\documentclass[a4paper,10pt,twoside]{scrartcl}
\usepackage{german,amssymb,amsmath,latexsym,epsfig,graphics,graphicx,mathrsfs,color,amsthm,cite}
\usepackage[latin1]{inputenc}
\usepackage{fancyhdr} 
\newtheorem{theorem}{Theorem}[section]

\newtheorem{definition}[theorem]{Definition}

\numberwithin{equation}{section}
\newcommand{\R}{{\mathbb R}}
\newcommand{\N}{{\mathbb N}}

\newcommand{\Zt}{{\mathscr{Z}^k([t_0,t_1])}}

\begin{document}
\thispagestyle{empty}
\rule{1.0\textwidth}{1.0pt}

\vspace*{5mm}
{\bf{\large Das Pontrjaginsche Maximumprinzip f"ur optimale Multiprozesse \\
            mit stetigen Zustandstrajektorien} \\[10mm]
Nico Tauchnitz} \\[25mm]
{\bf Vorwort} \\[2mm]
Die vorliegende Ausarbeitung ist eine Kurzfassung der Ergebnisse meiner Dissertation zu optimalen Multiprozessen
mit stetigen Zust"anden. In der Modellierung sind au"serdem Wechselkosten nicht vorgesehen.
Die Reduzierung auf stetige Zust"ande und die Vernachl"assigung von Wechselkosten schr"ankt die Aufgabenklasse im Rahmen
der Multiprozesse und hybriden Steuerungsprobleme ein.
Aber f"ur diese Aufgabenklasse gebe ich einen Vorschlag an,
wie man die Wechselstrategien zwischen den einzelnen Steuerungssystem auf elegante Weise als Steuervariable
in die Modellierung einbinden kann.
D.\,h. eine Vorgabe der Anzahl und der Zeitpunkte der Wechsel zwischen den verschiedenen Steuerungssystemen,
sowie die jeweilige Auswahl des Steuerungssystems ist nicht vorgeben.
Die Informationen "uber die optimale Wechselstrategie ergeben sich aus den Optimalit"atsbedingungen. \\[2mm]
Das Ergebnis f"ur diese Aufgabenklasse ist ein Pontrjaginsches Maximumprinzip f"ur starke lokale Optimalstellen.
Die notwendigen Bedingungen werden an illustrativen Beispielen demonstriert.
Dies ist ein Modell,
in dem zu jedem Zeitpunkt zwischen zwei verschiedenen Investitionsmodellen entschieden werden kann.
Unter besonderen Umst"anden f"uhrt dabei die Wechselstrategie zu einem Sattelpunkt.
Zus"atzlich wird ein zeitoptimales Modell gekoppelter Beh"alter untersucht.
Dieses Beispiel ist in meiner Dissertation noch nicht betrachtet wurden.
In diesem Beispiel k"onnen die Optimalit"atsbedingungen in der Weise entarten,
dass eine optimale Wechselstrategie nicht zu bestimmen ist.
Bei genauerem Blick auf das Modell kl"art sich dieser Fall sinnvoll auf. \\[5mm]
November 2015

\newpage
\lhead[\thepage \, Inhaltsverzeichnis]{Optimale Steuerung mit unendlichem Zeithorizont}
\rhead[Optimale Steuerung mit unendlichem Zeithorizont]{Inhaltsverzeichnis \thepage}
\tableofcontents

\newpage
\lhead[\thepage \, Einleitung]{Optimale Multiprozesse}
\rhead[Optimale Multiprozesse]{Einleitung \thepage}
\section{Einleitung} \label{SectionEinleitung}
Ein Multiprozess besteht aus einer gewissen endlichen Anzahl von einzelnen Steuerungssystemen,
welche sich jeweils aus individuellen Dynamiken, Zielfunktionalen und Steuerungsbereichen zusammensetzen.
Neben der Bestimmung der optimalen Steuerungen f"ur das jeweils gew"ahlte Steuerungssystem
ist das Ziel bei der Untersuchung von Multiprozessen die Anzahl der Wechsel
zwischen den einzelnen Steuerungssystemen, die konkrete Wahl des jeweiligen Steuerungssystems und diejenigen Zeitpunkte,
zu denen diese Wechsel stattfinden, zu bestimmen. \\
Als illustratives Beispiel f"uhren wir das Autofahren an:
Neben den kontinuierlichen Steuerungen ``beschleunigen'' und ``bremsen''
bildet die Wahl des konkreten Ganges einen wesentlichen Beitrag zur Minmierung des Treibstoffverbrauchs oder
zum Erreichen des Zielortes in k"urzester Zeit.
Die Auswahl des jeweiligen Ganges h"angt dabei z.\,B. von der momentanen Geschwindigkeit oder dem Streckenprofil oder
der Stra"sen- und Verkehrssituation ab.
Ma"sgeblich beeinflusst wird durch die Schaltvorg"ange das dynamische Beschleunigungsverhalten und der momentane Benzinverbrauch
oder die Geschwindigkeit. \\[2mm]
Eine anschauliche Beschreibung eines Multiprozesses wird durch die Auflistung der jeweilig ausgew"ahlten Steuerungssysteme gegeben.
Dazu zerlegt man das Zeitintervall $[0,T]$ zu den Wechselzeitpunkten in Teilintervalle,
$$0=t_0<t_1<...<t_N=T,$$
und beschreibt f"ur $j=1,...,N$ auf den Teilabschnitten die ausgew"ahlten Steuerungssysteme:
$$\int_{t_{j-1}}^{t_j} f_j\big(t,x(t),u(t)\big) \, dt, \quad
  \dot{x}(t) = \varphi_j\big(t,x(t),u(t)\big) \mbox{ f"ur } t \in (t_{j-1},t_j),\quad u(t) \in U \subseteq \R^{m_j}.$$
Der Multiprozess erh"alt (ohne die Beachtung weiterer Beschr"ankugen) dadurch die Gestalt
\begin{eqnarray*}
&& J\big(x(\cdot),u(\cdot)\big) = \sum_{j=1}^N \int_{t_{j-1}}^{t_j} f_j\big(t,x(t),u(t)\big) \, dt \to \inf, \\
&& \dot{x}(t) = \varphi_j\big(t,x(t),u(t)\big) \mbox{ f"ur } t \in (t_{j-1},t_j), \quad j=1,...,N, \\
&& u(t) \in U \subseteq \R^{m_j}, \quad j=1,...,N.
\end{eqnarray*}
An dieser Stelle muss kritisch hervorgehoben werden,
dass in dieser Darstellung eines Multiprozesses durch die Indizes $j=1,...,N$ die Struktur der Wechselstrategie vorgegeben ist.
D.\,h. die Anzahl der Wechsel und die jeweilige Auswahl des Steuerungssystems ist implizit bereits vorgegeben.
Die Beschreibung einer Wechselstrategie in Form einer ``echten'' Steuervariable wird in dieser Form nicht erreicht. \\[2mm]
Auf dieser Formulierung basiert die grundlegende Strategie zur Untersuchung von allgemeinen Multiprozessen
(insbesondere mit unstetigen Zust"anden und mit Wechselkosten):
{\em \glqq The maximum principle gives necessary conditions for [...] a solution [...].
The result only depends on comparing trajectories with the same switching strategy [...]\grqq} (Sussmann \cite{Sussmann}, S.\,330.).
Die Arbeit \cite{Sussmann} bezieht sich zwar auf allgemeinere Aufgaben als wir sie betrachten wollen,
aber die wesentlichen Beitr"age zu notwendigen Optimalit"atsbedingungen beruhen ebenso auf der Grundlage lediglich
gleichartige Wechselstrategien zu vergleichen
\cite{ClarkeVinter,ClarkeVinterII,Dmitruk,Galbraith,Garavello,Shaikh04,Shaikh07}. \\[2mm]
Wir zeigen in dieser Ausarbeitung f"ur Multiprozesse mit stetigen Zust"anden einen Weg auf,
die Wechselstrategien als Steuervariablen einzuf"uhren.
Damit sind wir in der Lage,
Multiprozesse auf starke lokale Optimalstellen zu untersuchen.
Das beinhaltet,
dass wir zur Bestimmung von notwendigen Optimalit"atsbedingungen Multiprozesse mit beliebig gearteten
Wechseltstrategien vergleichen. \\[2mm]
Wir beginnen diese Ausarbeitung mit der Einf"uhrung von geeigneten Zerlegungen des Zeitintervalls und der
Verkn"upfung dieser Variablen mit vektorwertigen Abbildungen.
Die Zerlegungen des Zeitintervalls bilden die Grundlage,
um im n"achsten Schritt die Wechselstrategien als Steuervariablen in die Beschreibung eines Multiprozesses einzuf"ugen,
sowie die konkrete Aufgabestellung und das Pontrjaginsche Maximumprinzip zu formulieren.
Die angegebenen notwendigen Bedingungen wenden wir anschlie"send auf ein einfaches Investitionsmodell an.
Danach weiten wir die Multiprozesse auf Aufgaben mit freien Endzeiten aus und betrachten abschlie"send die zeitoptimale Steuerung
eines gekoppelter Beh"alter. \\[2mm]
Wir gegeben f"ur das Pontrjaginsche Maximumprinzip nur Bemerkungen zum Beweis an.
Die ausf"uhrlichen Darlegungen zu optimalen Multiprozessen f"ur Aufgaben mit festen und freien Endzeiten sind in \cite{Tauchnitz}
zu finden.

\newpage
\lhead[\thepage \, Zerlegungen]{Optimale Multiprozesse}
\rhead[Optimale Multiprozesse]{Zerlegungen \thepage}
\section{$k$-fache Zerlegungen}
Es seien das Intervall $I \subseteq \R$ und $k \in \N$ gegeben.
\begin{definition}[$k$-fache Zerlegung] \label{DefinitionZerlegung} 
Unter einer $k$-fachen Zerlegung des Intervalls $I \subseteq \R$ verstehen wir ein endliches System
$\mathscr{A} = \{\mathscr{A}^1,...,\mathscr{A}^k\}$ von Lebesgue-me"sbaren Teilmengen von $I$ mit
$$\bigcup_{1\leq i\leq k} \mathcal{A}^{i} = I, \qquad
  \mathcal{A}^{i} \cap \mathcal{A}^{j} = \emptyset \mbox{ f"ur } i \not=j.$$
Es bezeichnet $\mathscr{Z}^k(I) = \{ \mathscr{A} \} = \big\{\{ \mathscr{A}^{i}\}_{1\leq i\leq k}\big\}$ die Menge aller
$k$-fachen Zerlegungen von $I$.
\end{definition}

Wir identifizieren die Elemente $\mathscr{A} \in \mathscr{Z}^k(I)$ durch die Vektorfunktion
$$\chi_{\mathscr{A}}(t) = \big(\chi_{\mathscr{A}^1}(t),...,\chi_{\mathscr{A}^k}(t) \big), \quad t \in I,$$
der charakteristischen Funktionen der Mengen ${\mathcal A}^{i}$.
Dann geh"ort jede Funktion $\chi_{\mathscr{A}}(\cdot)$ zu der Menge
$$\mathscr{Y}^k(I) = \bigg\{ y(\cdot) \in L_\infty(I,\R^k)\,\bigg|\, y(t)=\big(y^1(t),...,y^k(t)\big), y^{i}(t) \in \{0,1\},
                             \sum_{i=1}^ky^{i}(t)=1\bigg\}.$$
Die Menge $\mathscr{Z}^k(I)$ der $k$-fachen Zerlegungen bzw. $\mathscr{Y}^k(I)$ der charakteristischen Vektorfunktionen
repr"asentieren s"amtliche Wechselstrategien in Multiprozessen mit $k$ Steuerungssystemen. \\[2mm]
Wir betrachten nun $k$ Funktion $h^{i}: I \times \R^n \times \R^{m_i} \to \R^s$, und fassen diese zur
Vektorfunktion $h=(h^1,...,h^k)$ zusammen.
Die Verkn"upfung mit einer $k$-fachen Zerlegung legen wir wie folgt fest:
$$\mathscr{A} \circ h(t,x,u) = \chi_{\mathscr{A}}(t) \circ h(t,x,u)
                             = \sum_{i=1}^k \chi_{\mathscr{A}^{i}}(t) \cdot h^{i}(t,x,u^{i}),
  \quad u=(u^1,...,u^k) \in \R^{m_1+...+m_k}.$$
Sind die Funktionen $h^{i}$ differenzierbar, dann setzen wir f"ur die partielle Ableitung bzgl. $t$ und $x$:
\begin{eqnarray*}
\mathscr{A} \circ h_t(t,x,u) &=& \chi_{\mathscr{A}}(t) \circ h_t(t,x,u)
                             = \sum_{i=1}^k \chi_{\mathscr{A}^{i}}(t) \cdot h^{i}_t(t,x,u^{i}), \\
\mathscr{A} \circ h_x(t,x,u) &=& \chi_{\mathscr{A}}(t) \circ h_x(t,x,u)
                             = \sum_{i=1}^k \chi_{\mathscr{A}^{i}}(t) \cdot h^{i}_x(t,x,u^{i}).
\end{eqnarray*}
Die $k$-fachen Zerlegungen besitzen grundlegende Eigenschaften,
von denen wir einige angeben werden 
F"ur $\chi_{\mathscr{A}}(\cdot),\chi_{\mathscr{B}}(\cdot) \in \mathscr{Y}^k(I)$,
$g,h:I \times \R^n \times \R^m \to \R^s$ und $\alpha \in \R$ gelten offenbar:
\begin{itemize}
\item[(a)] $\chi_{\mathscr{A}}(t) \circ [\alpha \cdot h(t,x,u)] = \alpha \cdot [\chi_{\mathscr{A}}(t) \circ h(t,x,u)]$,
\item[(b)] $\chi_{\mathscr{A}}(t) \circ [g(t,x,u) + h(t,x,u)]
             = [\chi_{\mathscr{A}}(t) \circ g(t,x,u)]+[\chi_{\mathscr{A}}(t) \circ h(t,x,u)]$,
\item[(c)] $[\chi_{\mathscr{A}}(t)+ \chi_{\mathscr{B}}(t)] \circ h(t,x,u)
             = [\chi_{\mathscr{A}}(t) \circ h(t,x,u)]+ [\chi_{\mathscr{B}}(t) \circ h(t,x,u)]$.
\end{itemize}
Eine weitere Eigenschaft, die in Bezug zu Nadelvariationen von Wechselstrategien steht, ist:
\begin{itemize}
\item[(d)] Es seien $\chi_{\mathscr{A}}(\cdot), \chi_{\mathscr{B}}(\cdot) \in \mathscr{Y}^k(I)$ und $M \subseteq I$ me"sbar.
           Wir setzen
           $$y(t)= \chi_{\mathscr{A}}(t) + \chi_M(t) \big(\chi_{\mathscr{A}}(t)-\chi_{\mathscr{B}}(t)\big).$$
           Dann gilt $y(\cdot) \in \mathscr{Y}^k(I)$.
\end{itemize}

\newpage
\lhead[\thepage \, Aufgabenstellung]{Optimale Multiprozesse}
\rhead[Optimale Multiprozesse]{Aufgabenstellung \thepage}
\section{Die Aufgabenstellung}  \label{KapitelHybridFest}
Es sei $k \in \N$ die Anzahl der verschiedenen gegebenen Steuerungssysteme,
die individuelle Integranden $f^{i}(t,x,u^{i}): \R \times \R^n \times \R^{m_i} \to \R$,
Dynamiken $\varphi^{i}(t,x,u^{i}): \R \times \R^n \times \R^{m_i} \to \R^n$
und Steuerbereiche $U^{i} \subseteq \R^{m_i}$ besitzen k"onnen.
Ferner seien
$$g_j(t,x):\R \times \R^n \to \R, \quad j=1,...,l, \qquad h_\iota(x):\R^n \to \R^{s_\iota }, \quad \iota = 0,1.$$
Wir fassen die Steuerungssysteme durch die Setzungen
$$f=(f^1,...,f^k),\qquad \varphi=(\varphi^1,...,\varphi^k), \qquad U=U^1 \times ... \times U^k$$
zusammen.
Mit den Vorbetrachtungen des letzten Abschnitts formulieren wir "uber dem Zeitintervall $[t_0,t_1]$
die Aufgabe des optimalen Multiprozesses mit Zustandsbeschr"ankungen wie folgt:
\begin{eqnarray}
&& \label{HybrideAufgabe1} J\big(x(\cdot),u(\cdot),\mathscr{A}\big)
           = \int_{t_0}^{t_1} \chi_{\mathscr{A}}(t) \circ f\big(t,x(t),u(t)\big) dt \to \inf, \\
&& \label{HybrideAufgabe2} \dot{x}(t) = \chi_{\mathscr{A}}(t) \circ \varphi\big(t,x(t),u(t) \big), \\
&& \label{HybrideAufgabe3} g_j\big(t,x(t)\big) \leq 0, \quad t \in [t_0,t_1], \quad j=1,...,l, \\
&& \label{HybrideAufgabe4} h_0\big( x(t_0) \big) = 0, \quad h_1\big( x(t_1) \big) = 0, \\
&& \label{HybrideAufgabe5} u(t) \in U= U^1 \times ... \times U^k, \quad U^{i} \not= \emptyset, \quad \mathscr{A} \in \Zt.
\end{eqnarray}
Wir treffen folgende Annahmen in der Aufgabe (\ref{HybrideAufgabe1})--(\ref{HybrideAufgabe5}): 
\begin{itemize}
\item[{\bf (A$_1$)}] F"ur $i=1,...,k$ sind die Abbildungen $f^{i}(t,x,u^{i})$, $\varphi^{i}(t,x,u^{i})$ stetig in der Gesamtheit
      der Variablen, stetig differenzierbar bzgl. $x$ und es seien die Abbildungen $f^{i}_x(t,x,u^{i}), \varphi^{i}_x(t,x,u^{i})$
      stetig in der Gesamtheit der Variablen.
\item[{\bf (A$_2$)}] S"amtliche Funktionen $g_j(t,x)$, $h_\iota(x)$ sind stetig, stetig nach $x$ differenzierbar und es seien die 
      Abbildungen $g_{jx}(t,x)$ stetig in der Gesamtheit der Variablen.
\end{itemize}
Wir nennen ein Tripel $\big(x(\cdot),u(\cdot),\mathscr{A}\big)$ mit
$x(\cdot) \in W^1_\infty\big([t_0,t_1],\R^n\big)$, $u(\cdot) \in L_\infty\big([t_0,t_1],U\big)$ und Wechselstrategie
$\mathscr{A} \in \Zt$ einen Multiprozess.
Ein Multiprozess $\big(x(\cdot),u(\cdot),\mathscr{A}\big)$ hei"st zul"assig in der
Aufgabe (\ref{HybrideAufgabe1})--(\ref{HybrideAufgabe5}),
wenn auf dem Intervall $[t_0,t_1]$ die Funktion $x(\cdot)$ fast "uberall der Gleichung
$$\dot{x}(t) = \chi_{\mathscr{A}}(t) \circ \varphi \big(t,x(t),u(t)\big)
             = \sum_{i=1}^k \chi_{\mathscr{A}^{i}}(t) \cdot \varphi^{i}\big(t,x(t),u^{i}(t)\big)$$
gen"ugt, f"ur alle $t \in [t_0,t_1]$ die Zustandsbeschr"ankungen 
$g_j\big(t,x(t)\big) \leq 0$, $j=1,...,l,$
erf"ullt sind und in den Endpunkten $t_0,t_1$ die Randbedingungen
$h_0\big( x(t_0) \big) = 0$, $h_1\big( x(t_1) \big) = 0$ gelten. \\
Einen zul"assigen Multiprozess $\big(x_*(\cdot),u_*(\cdot),\mathscr{A}_*\big)$ nennen wir ein starkes lokales Minimum,
wenn eine Zahl $\varepsilon > 0$ derart existiert,
dass f"ur jeden zul"assigen Multiprozess $\big(x(\cdot),u(\cdot),\mathscr{A}\big)$ mit
$$\| x(\cdot)-x_*(\cdot) \|_\infty < \varepsilon$$
folgende Ungleichung gilt:
$$J\big(x(\cdot),u(\cdot),\mathscr{A}\big) \geq J\big(x_*(\cdot),u_*(\cdot),\mathscr{A}_*\big).$$

\newpage
\lhead[\thepage \, Maximumprinzip]{Optimale Multiprozesse}
\rhead[Optimale Multiprozesse]{Maximumprinzip \thepage}
\section{Das Pontrjaginsche Maximumprinzip}
In der Aufgabe (\ref{HybrideAufgabe1})--(\ref{HybrideAufgabe5})
definieren wir f"ur die einzelnen Steuerungssysteme die partiellen Pontrjaginschen Funktionen
$H^{i}:\R \times \R^n \times \R^{m_i} \times \R^n \times \R \to \R$
standardgem"a"s durch
$$H^{i}(t,x,u^{i},p,\lambda_0) = \langle p , \varphi^{i}(t,x,u^{i}) \rangle - \lambda_0 f^{i}(t,x,u^{i}).$$
Dann hei"st f"ur die Aufgabe optimaler Multiprozesse die Funktion
\begin{eqnarray*}
    \chi_{\mathscr{A}}(t) \circ H(t,x,u,p,\lambda_0\big)
&=& \sum_{i=1}^k \chi_{\mathscr{A}^{i}}(t) \cdot H^{i}(t,x,u^{i},p,\lambda_0) \\
&=& \sum_{i=1}^k \chi_{\mathscr{A}^{i}}(t) \cdot \big( \langle p , \varphi^{i}(t,x,u^{i}) \rangle - \lambda_0 f^{i}(t,x,u^{i}) \big)
\end{eqnarray*}
mit $H:\R \times \R^n \times \R^{m_1 + ... + m_k} \times \R^n \times \R \to \R^k$ die Pontrjaginsche Funktion.

\begin{theorem} [Pontrjaginsches Maximumprinzip] \label{SatzPMPhybridFest}
Es seien die Annahmen {\bf (A$_1$)}, {\bf (A$_2$)} erf"ullt. \\[1mm]
Ist $\big(x_*(\cdot),u_*(\cdot),\mathscr{A}_*\big)$ ein starkes lokales Minimum in der Aufgabe
(\ref{HybrideAufgabe1})--(\ref{HybrideAufgabe5}),
dann existieren eine Zahl $\lambda_0 \geq 0$, Vektoren $l_0 \in \R^{s_0},\, l_1 \in \R^{s_1}$,
eine Vektorfunktion $p(\cdot)$, $p: [t_0,t_1] \to \R^n,$
und auf den Mengen $T_j = \big\{ t \in [t_0,t_1] \,\big|\, g_j\big( t,x_*(t)\big) = 0 \big\}$
konzentrierte nichtnegative regul"are Borelsche Ma"se $\mu_j$, $j=1,...,l,$ (wobei diese Gr"o"sen nicht gleichzeitig verschwinden) derart,
dass
\begin{enumerate}
\item[a)] die Vektorfunktion $p(\cdot)$ der Integralgleichung
          \begin{eqnarray} 
          p(t) &=& - {h_1'}^T\big(x_*(t_1)\big) \, l_1
              + \int_t^{t_1} \chi_{\mathscr{A}_*}(\tau) \circ H_x\big(\tau,x_*(\tau),u_*(\tau),p(\tau),\lambda_0\big) d\tau \nonumber \\
               & & \label{SatzHybrideAufgabePMP1}  - \sum_{j=1}^l \int_t^{t_1} g_{jx}\big( \tau,x_*(\tau)\big) d\mu_j(\tau)
          \end{eqnarray}
          gen"ugt und die Transversalit"atsbedingung
          \begin{equation} \label{SatzHybrideAufgabePMP2} 
          p(t_0) = {h_0'}^T\big(x_*(t_0)\big) l_0
          \end{equation}
          erf"ullt;
\item[b)] f"ur fast alle $t$ aus $[t_0,t_1]$ die Maximumbedingung
          \begin{equation} \label{SatzHybrideAufgabePMP3} 
          \chi_{\mathscr{A}_*}(t) \circ H\big(t,x_*(t),u_*(t),p(t),\lambda_0\big)
          = \max_{\substack{u^{i} \in U^{i}\\ 1\leq i\leq k}} H^{i}\big(t,x_*(t),u^{i},p(t),\lambda_0\big).
          \end{equation}
          gilt.
\end{enumerate}
\end{theorem}

In den Aufgaben ohne Zustandsbeschr"ankungen ist $p(\cdot)$ eine absolutstetige Funktion.
Liegen aber Phasenbeschr"ankungen vor, so kann die adjungierte Funktion $p(\cdot)$ Unstetigkeiten aufweisen,
da in Gleichung (\ref{SatzHybrideAufgabePMP1}) Integrale bzgl. der Ma"se $\mu_j$ vorkommen.
Sie ist jedoch stets eine Funktion beschr"ankter Variation, die linksseitig stetig ist
(vgl. \cite{Ioffe}, S.\,209). \\[2mm]
{\bf Beweisskizze} Wir halten uns an die Beweismethode von Ioffe\& Tichomirov f"ur regul"ar lokalkonvexe Aufgaben
(\cite{Ioffe}, S.\,200--224).
Grundlage bildet dabei Satz 1, S.\,201. \\
Die Voraussetzungen a) und c) von Satz 1 folgen ebenso wie in \cite{Ioffe}, 
da jedes einzelne Steuerungssystem die klassischen Voraussetzungen an Optimalsteuerungsprobleme erf"ullt und zudem die Relation
$$\|\chi_{\mathscr{A}}(t) \circ h(t,x,u)\| \leq \sum_{i=1}^k \|h^{i}(t,x,u^{i})\|$$
f"ur die individuellen Abbildungen und jede Wechselstrategie gilt. \\
Es bleibt die Existenz der Abbildung unter Voraussetzung b) von Satz 1 nachzuweisen.
Die Grundlage daf"ur liefern die Mengenfamilien $M_1(\alpha),..., M_{\kappa_0}(\alpha)$ mit den Eigenschaften
(vgl. \cite{Ioffe}, Lemma\,2, S.\,217)
$$|M_\kappa(\alpha)| = \alpha (t_1-t_0), \quad M_\kappa(\alpha') \subseteq M_\kappa(\alpha),
   \quad M_\kappa(\alpha) \cap M_{\kappa'}(\alpha') = \emptyset.$$
Wegen der Disjunktheit der Mengenfamilien f"ur verschiedene Indizes
k"onnen wir zu gegebenen $u_\kappa(\cdot) \in L_\infty([t_0,t_1],U)$ die Steuerung $u_\alpha(\cdot) \in L_\infty([t_0,t_1],U)$ gem"a"s
$$u_\alpha(t)=u_*(t) + \sum_{\kappa=1}^{\kappa_0} \chi_{M_\kappa(\alpha_\kappa)}(t) \cdot \big( u_\kappa(t)-u_*(t) \big)$$
definieren.
Diese Nadelvariation der Steuerung $u_*(\cdot)$ bzgl. me"sbarer Mengen in Richtung der Steuerung $u_\kappa(\cdot)$
besitzt f"ur jede wohldefinierte $h=h(t,x,u)$ die Eigenschaft
$$h\big(t,x,u_\alpha(t) \big)=h \big(t,x,u_*(t) \big)+
  \sum_{\kappa=1}^{\kappa_0} \chi_{M_\kappa(\alpha_\kappa)}(t) \cdot \Big( h\big(t,x,u_\kappa(t)\big)- h\big(t,x,u_*(t) \big) \Big).$$
D.\,h., die Nadelvariationen bzgl. der impliziten Steuervariable wirken sich entsprechend als Nadelvariation
auf den Wert der Abbildung $h$ aus. \\
Mit der Einf"uhrung der Wechselstrategien als zus"atzliche Steuervariable ist diese Eigenschaft f"ur Multiprozesse sicherzustellen,
damit Voraussetzung b) von Satz 1 in \cite{Ioffe} erf"ullt ist.
Auf Basis der Eigenschaften, dass f"ur $\chi_{\mathscr{A}}(\cdot), \chi_{\mathscr{B}}(\cdot) \in \mathscr{Y}^k([t_0,t_1])$
und f"ur me"sbare Teilmengen $M \subseteq [t_0,t_1]$ die Funktion
$$y(t)= \chi_{\mathscr{A}}(t) + \chi_M(t) \big(\chi_{\mathscr{A}}(t)-\chi_{\mathscr{B}}(t)\big)$$
wieder eine Wechselstrategie beschreibt, l"asst sich die Beziehung
\begin{eqnarray}
\lefteqn{\chi_{\mathcal{A}_\alpha}(t) \circ h \big(t,x,u_\alpha(t) \big)
    - \chi_{\mathcal{A}_*}(t) \circ h \big(t,x,u_*(t) \big)} \nonumber \\
&=& \label{EigenschaftNadelvariationHybrid1}
    \sum_{\kappa=1}^{\kappa_0} \chi_{M_\kappa(\alpha_\kappa)}(t) \cdot \Big(  \chi_{\mathcal{A}_\kappa}(t) \circ h \big(t,x,u_\kappa(t) \big) 
    - \chi_{\mathcal{A}_*}(t) \circ h \big(t,x,u_*(t) \big) \Big)
\end{eqnarray}
und damit das Maximumprinzip \ref {SatzPMPhybridFest} (ausf"uhrlich in \cite{Tauchnitz}) zeigen. \hfill $\blacksquare$ 

\newpage
\lhead[\thepage \, Investitionsmodell]{Optimale Multiprozesse}
\rhead[Optimale Multiprozesse]{Investitionsmodell \thepage}
\section{Ein Investitionsmodell} \label{AbschnittBeispielhybrid}
Wir untersuchen ein Investitionsmodell,
dass sich aus den \glqq consumption vs. investement\grqq -Modellen in \cite{Seierstad},
S.\,82 und S.\,130, zusammensetzt:
\begin{eqnarray}
&& \label{Beispielhybrid1} J\big(x(\cdot),u(\cdot),\mathscr{A} \big)
                            = \int_0^T \chi_{\mathscr{A}}(t) \circ f\big(t,x(t),u(t)\big) dt \rightarrow \sup, \\
&& \label{Beispielhybrid2} \dot{x}(t) = \chi_{\mathscr{A}}(t) \circ \varphi\big(t,x(t),u(t)\big), \quad x(0)=x_0>0, \\
&& \label{Beispielhybrid3} u(t) \in [0,1], \quad \mathscr{A}=\{ \mathscr{A}^1,\mathscr{A}^2 \} \in \mathscr{Z}^2([0,1]),
\end{eqnarray}
mit den einzelnen Steuerungssystemen
$$\begin{array}{ll}
     f^1(t,x,u^1) = (1-u^1)x,        & \varphi^1(t,x,u^1) = u^1 x, \\[1mm]
     f^2(t,x,u^2) = (1-u^2)x^\alpha, & \varphi^2(t,x,u^2) = u^2 x^\alpha,
  \end{array} \quad U^1=U^2=[0,1]$$
und mit den Modellparametern
\begin{equation} \label{Beispielhybrid4}
\alpha \in (0,1) \mbox{ konstant}, \qquad T \mbox{ fest mit } T > \max \left\{1, \frac{x_0^{1-\alpha}}{\alpha} \right\}.
\end{equation}
Wir wenden Theorem \ref{SatzPMPhybridFest} an:
Angenommen, es ist $\lambda_0=0$.
Dann folgt $p(t),q(t)\equiv 0$, im Widerspruch zu Satz \ref{SatzPMPhybridFest}.
Daher sei $\lambda_0=1$.
Dann gilt nach (\ref{SatzHybrideAufgabePMP3}) f"ur fast alle $t \in [0,T]$:
\begin{eqnarray*}
\lefteqn{\chi_{\mathscr{A}_*}(t) \circ H\big(t,x_*(t),u_*(t),p(t),1\big)} \\
&=& \max_{i \in \{1,2\}} \Big\{ \max_{u^1 \in [0,1]} \big[ \big( p(t) -1 \big) u^1 + 1 \big] \cdot x_*(t),
                \max_{u^2 \in [0,1]} \big[ \big( p(t) -1 \big) u^2 + 1 \big] \cdot x_*^\alpha(t) \Big\}.
\end{eqnarray*}
Dies k"onnen wir weiterhin in die Form
\begin{equation}
\label{Beispielhybrid5} \max_{v \in [0,1]} \big[ \big( p(t) -1 \big) v \big] \cdot \max_{i \in \{1,2\}}\{ x_*(t) , x_*^\alpha(t) \}
\end{equation}
bringen.
Daraus erhalten wir f"ur die optimale Investitionsrate und Wechselstrategie
$$u^1_*(t)=u^2_*(t)=v_*(t) = \left\{ \begin{array}{ll}
              1, & p(t) < 1, \\
              0, & p(t) > 1, \end{array} \right. \quad
  \chi_{\mathscr{A}_*}(t) = \left\{ \begin{array}{ll}
               ( 1,0 ), & x_*(t) > 1, \\
               ( 0,1 ), & 0 < x_*(t) < 1. \end{array} \right.$$
Die Funktion $p(\cdot)$ ist die L"osung der adjungierten Gleichung (\ref{SatzHybrideAufgabePMP1})
\begin{equation} \label{Beispielhybrid6}
\dot{p}(t) = \left\{ \begin{array}{ll} - \big[ \big( p(t) -1 \big) v_*(t) + 1 \big],    & t \in \mathscr{A}_1, \\[1mm]
          - \big[ \big( p(t) -1 \big) v_*(t) + 1 \big] \cdot \alpha x_*^{\alpha-1}(t), & t \in \mathscr{A}_2, \end{array} \right.
\quad p(T)=0.
\end{equation}
Betrachten wir die Steuerungssysteme f"ur sich,
so sind durch $t_1=T-1$ bzw. $t_2=\alpha T -x_0^{1-\alpha}$ die Zeitpunkte f"ur den optimalen Wechsel von vollst"andiger
Investition in komplette Kosumption gegeben.
Au"serdem ist f"ur die L"osung $x(\cdot)$ der Differentiagleichung $\dot{x}(t)=x^\alpha(t)$ mit $x_0 \in (0,1)$ der
Zeitpunkt $\sigma$, in dem $x(\sigma)=1$ gilt, durch $\sigma=\frac{1-x_0^{1-\alpha}}{1-\alpha}$ bestimmt.
Im Weiteren seien also
$$t_1=T-1,  \qquad t_2=\alpha T -x_0^{1-\alpha}, \qquad \sigma=\frac{1-x_0^{1-\alpha}}{1-\alpha}.$$
\begin{enumerate}
\item[a)] Sei $x_0 \geq 1$: In diesem Fall lautet die L"osung
          \begin{eqnarray*}
          && x_*(t) = \left\{ \begin{array}{ll} x_0 \cdot e^t, & t \in [0,t_1), \\
                                                x_0 \cdot e^{t_1}, & t \in [t_1,T], \end{array} \right.
             \qquad u^1_*(t) = \left\{ \begin{array}{ll} 1, & t \in [0,t_1), \\
                                                        0, & t \in [t_1,T], \end{array} \right. \\
          && \chi_{\mathscr{A}_*}(t) = (1,0), \qquad J\big(x_*(\cdot),u_*(\cdot),\mathscr{A}_*\big) = x_0 \cdot e^{t_1}.          
          \end{eqnarray*}
\item[b)] Seien $x_0<1$ und $T < \frac{1 - \alpha x_0^{1-\alpha}}{\alpha (1-\alpha)}$:
          Dann lautet die L"osung
          \begin{eqnarray*}
          && y_*(t) = \left\{ \begin{array}{ll}
                      \big[ (1-\alpha)t +  x_0^{1-\alpha} \big]^\frac{1}{1-\alpha}, & t \in [0,t_2), \\
                      \big[ \alpha(1-\alpha)T + \alpha x_0^{1-\alpha} \big]^\frac{1}{1-\alpha}, & t \in [t_2,T],
                      \end{array} \right.
             \quad v^2_*(t) = \left\{ \begin{array}{ll} 1, & t \in [0,t_2), \\
                                                        0, & t \in [t_2,T], \end{array} \right. \\
          && \chi_{\mathscr{B}_*}(t) = (0,1), \qquad 
             J\big(y_*(\cdot),v_*(\cdot),\mathscr{B}_*\big) = 
             \alpha^\frac{\alpha}{1-\alpha} \cdot \big[ (1-\alpha)T + x_0^{1-\alpha} \big]^\frac{1}{1-\alpha}.          
          \end{eqnarray*}
\item[c)] Seien $x_0<1$ und $T > \frac{1 - \alpha x_0^{1-\alpha}}{\alpha (1-\alpha)}$:
          Nun erhalten wir die L"osung
          \begin{eqnarray*}
          && z_*(t) = \left\{ \begin{array}{ll}
                      \big[ (1-\alpha)t +  x_0^{1-\alpha} \big]^\frac{1}{1-\alpha}, & t \in [0,\sigma), \\
                      e^{t-\sigma}, & t \in [\sigma,t_1), \\
                      e^{t_1-\sigma}, & t \in [t_1,T], \end{array} \right.                      
             \quad \begin{array}{ll}
                   & w^2_*(t)= \;\; 1, \; t \in [0,\sigma), \\
                   & w^1_*(t)= \left\{ \begin{array}{ll} 1, & t \in [\sigma,t_1), \\
                                                        0, & t \in [t_1,T], \end{array} \right.
                   \end{array} \\
          && \chi_{\mathscr{C}_*}(t) = \left\{ \begin{array}{ll} (0,1), & t \in [0,\sigma), \\
                                                        (1,0), & t \in [\sigma,T], \end{array} \right.
             \qquad J\big(z_*(\cdot),w_*(\cdot),\mathscr{C}_*\big) = e^{T-1-\sigma}.        
          \end{eqnarray*}
\item[d)] Seien $x_0<1$ und $T = \frac{1 - \alpha x_0^{1-\alpha}}{\alpha (1-\alpha)}$:
          Wegen $\sigma=t_2$ erf"ullen beide Multiprozesse $\big(y_*(\cdot),v_*(\cdot),\mathscr{B}_*\big)$ und
          $\big(z_*(\cdot),w_*(\cdot),\mathscr{C}_*\big)$ alle Bedingungen von Theorem \ref{SatzPMPhybridFest}. \\
          Sei $n_0 \in \N$ mit $n_0 >\alpha/(1-\alpha)$.
          F"ur $n \geq n_0$ betrachten wir die Steuerungen
          $$v_n(t)= v_*(t) + \chi_{[\sigma,\sigma+\frac{1}{n})}(t) \cdot \big(w_*(t)-v_*(t)\big),
            \quad \chi_{\mathscr{B}_n}(t)= \chi_{\mathscr{C}_*}(t),$$
          die wie im Fall c) einen Wechsel des Steuerungssystems zum Zeitpunkt $t=\sigma=t_2$ und eine verl"angerte Investitionsphase
          vorgeben.
          F"ur den zugeh"origen Kapitalbestand $y_n(\cdot)$ gilt $\|y_n(\cdot)-y_*(\cdot)\|_\infty=e^{\frac{1}{n}}-1$
          und f"ur den Wert des Zielfunktionals
          \begin{eqnarray*}
          \lefteqn{J\big(y_n(\cdot),v_n(\cdot),\mathscr{B}_n\big)-J\big(y_*(\cdot),v_*(\cdot),\mathscr{B}_*\big)
           =e^{\frac{1}{n}}\bigg(T-\sigma-\frac{1}{n}\bigg)- (T-\sigma)} \\
          && >\bigg(1+\frac{1}{n}\bigg)\bigg(\frac{1}{\alpha}-\frac{1}{n}\bigg)-\frac{1}{\alpha}
             =\frac{1}{n}\bigg(\frac{1}{\alpha}-1-\frac{1}{n}\bigg)>0 \quad\mbox{ f"ur alle } n \geq n_0.
          \end{eqnarray*}
          Daher stellt $\big(y_*(\cdot),v_*(\cdot),\mathscr{B}_*\big)$ in diesem Fall kein starkes lokales Optimum dar.
          Der Multiprozess $\big(z_*(\cdot),w_*(\cdot),\mathscr{C}_*\big)$ ist die eindeutige L"osung.
\end{enumerate}

\newpage
\lhead[\thepage \, Freie Endzeiten]{Optimale Multiprozesse}
\rhead[Optimale Multiprozesse]{Freie Endzeiten \thepage}
\section{Multiprozesse mit freiem Anfangs- und Endzeitpunkt} \label{KapitelFreieZeit}
Wir betrachten in diesem Abschnitt Multiprozesse mit freiem Anfangs- und Endzeitpunkt.
Im Gegensatz zur Aufgabe (\ref{HybrideAufgabe1})--(\ref{HybrideAufgabe5}) bedeutet dies,
dass die Abbildungen $h_0$ bzw. $h_1$, die die Start- und Zielmannigfaltigkeit definieren, zeitabh"angig sind:
$$h_\iota (t,x) : \R \times \R^n \to \R^{s_\iota }, \quad \iota = 0,1.$$
Zur Behandlung dieser Aufgaben gehen wir analog zu \cite{Ioffe}, S.\,210--213,
vor und "uberf"uhren das Problem mit freier Zeit mittels einer Substitution der Zeit in eine Aufgabe auf festem Zeitintervall.
Diese Substitution der Zeit macht es unumg"anglich,
dass die eingehenden Funktionen und Abbildungen bzgl. der Zeitvariablen $t$ differenzierbar sind. \\[2mm]
Wir "ubernehmen alle in Abschnitt \ref{KapitelHybridFest} getroffenen Bezeichnungen und Festlegungen, und betrachten im Folgenden
Multiprozesse mit freiem Zeitintervall $[t_0,t_1]$:
\begin{eqnarray}
&& \label{HybrideAufgabeFrei1} J\big(x(\cdot),u(\cdot),\mathscr{A} \big)
           = \int_{t_0}^{t_1} \chi_{\mathscr{A}}(t) \circ f\big(t,x(t),u(t)\big) dt \to \inf, \\
&& \label{HybrideAufgabeFrei2} \dot{x}(t) = \chi_{\mathscr{A}}(t) \circ \varphi \big(t,x(t),u(t) \big), \\
&& \label{HybrideAufgabeFrei3} g_j\big(t,x(t)\big) \leq 0, \quad t \in [t_0,t_1], \quad j=1,...,l, \\
&& \label{HybrideAufgabeFrei4} h_0\big( t_0,x(t_0) \big) = 0, \quad h_1\big( t_1,x(t_1) \big) = 0, \\
&& \label{HybrideAufgabeFrei5} u(t) \in U= U^1 \times ... \times U^k, \quad U^{i} \not= \emptyset, \quad \mathscr{A} \in \Zt.
\end{eqnarray}
Wir treffen folgende Annahmen in der Aufgabe (\ref{HybrideAufgabeFrei1})--(\ref{HybrideAufgabeFrei5}):
\begin{itemize}
\item[{\bf (B$_1$)}] F"ur $i=1,...,k$ sind die Abbildungen $f^{i}(t,x,u^{i})$, $\varphi^{i}(t,x,u^{i})$ stetig in der Gesamtheit
      der Variablen, stetig differenzierbar bzgl. $t$ und $x$ und es seien alle partiellen Ableitungen stetig
      in der Gesamtheit der Variablen.
\item[{\bf (B$_2$)}] S"amtliche Funktionen $g_j(t,x)$, $h_\iota(t,x)$ sind stetig und stetig differenzierbar.
\end{itemize}
Wir nennen ein Quadrupel $\big(x(\cdot),u(\cdot),\mathscr{A},[t_0,t_1]\big)$ mit
$x(\cdot) \in W^1_\infty\big([t_0,t_1],\R^n\big)$, $u(\cdot) \in L^\infty\big([t_0,t_1],U\big)$, $\mathscr{A} \in \Zt$
und $[t_0,t_1] \subset \R$ einen Multiprozess mit freier Zeit.
Ein Multiprozess mit freier Zeit hei"st zul"assig in der Aufgabe (\ref{HybrideAufgabeFrei1})--(\ref{HybrideAufgabeFrei5}),
wenn auf dem Intervall $[t_0,t_1]$ die Funktion $x(\cdot)$ fast "uberall der Gleichung
$$\dot{x}(t) = \chi_{\mathscr{A}}(t) \circ \varphi \big(t,x(t),u(t) \big)$$
gen"ugt, f"ur alle $t \in [t_0,t_1]$ die Zustandsbeschr"ankungen $g_j\big(t,x(t)\big) \leq 0$, $j=1,...,l,$
erf"ullt sind und in den Endpunkten $t_0,t_1$ die Randbedingungen
$h_0\big( t_0,x(t_0) \big) = 0$, $h_1\big( t_1,x(t_1) \big) = 0$ gelten. \\
Ein zul"assiger Multiprozess $\big(x_*(\cdot),u_*(\cdot),\mathscr{A}_*,[t_{0*},t_{1*}]\big)$ ist ein starkes lokales Minimum,
wenn eine Zahl $\varepsilon > 0$ derart existiert,
dass f"ur jeden anderen zul"assigen Steuerprozess $\big(x(\cdot),u(\cdot),\mathscr{A},[t_0,t_1]\big)$ mit den Eigenschaften
$$|t_0 - t_{0*}| < \varepsilon, \quad |t_1 - t_{1*}| < \varepsilon$$
und
$$\| x(\cdot)-x_*(\cdot) \|_\infty < \varepsilon \quad\mbox{ f"ur jedes } t \in [t_{0*},t_{1*}] \cap [t_0,t_1]$$
folgende Ungleichung gilt:
$$J\big(x(\cdot),u(\cdot),\mathscr{A}\big) \geq J\big(x_*(\cdot),u_*(\cdot),\mathscr{A}_*\big).$$

In der Aufgabe (\ref{HybrideAufgabeFrei1})--(\ref{HybrideAufgabeFrei5})
bezeichnen
\begin{eqnarray*}
    \chi_{\mathscr{A}}(t) \circ H(t,x,u,p,\lambda_0\big)
&=& \sum_{i=1}^k \chi_{\mathscr{A}^{i}}(t) \cdot H^{i}(t,x,u^{i},p,\lambda_0) \\
    \mathscr{H}(t,x,p,\lambda_0)
&=& \max_{\substack{u^{i} \in U_i\\ 1\leq i\leq k}} H_i(t,x,u^{i},p,\lambda_0)
\end{eqnarray*}
die Pontrjagin-Funktion bzw. die Hamilton-Funktion.

\begin{theorem} [Pontrjaginsches Maximumprinzip] \label{SatzFreieZeitPMP}
Es seien die Annahmen {\bf (B$_1$)}, {\bf (B$_2$)} erf"ullt. \\[1mm]
Ist $\big(x_*(\cdot),u_*(\cdot),\mathscr{A}_*,[t_{0*},t_{1*}]\big)$ ein starkes lokales Minimum in der Aufgabe
(\ref{HybrideAufgabeFrei1})--(\ref{HybrideAufgabeFrei5}),
dann existieren eine Zahl $\lambda_0 \geq 0$, Vektoren $l_0 \in \R^{s_0},\, l_1 \in \R^{s_1}$, eine Vektorfunktion
$p(\cdot), p: [t_0,t_1] \to \R^n,$ und auf den Mengen $T_j = \big\{ t \in [t_0,t_1] \,\big|\, g_j\big( t,x_*(t)\big) = 0 \big\}$
konzentrierte nichtnegative regul"are Borelsche Ma"se $\mu_j$, $j=1,...,l,$ (wobei diese Gr"o"sen nicht gleichzeitig verschwinden) derart,
dass
\begin{enumerate}
\item[a)] die Vektorfunktion $p(\cdot)$ der Integralgleichung
          \begin{eqnarray} 
          p(t) &=& -{h_1}_x^T\big(t_{1*},x_*(t_{1*})\big) l_1 +
               \int_t^{t_{1*}} \chi_{\mathscr{A}_*}(\tau) \circ H_x\big(\tau,x_*(\tau),u_*(\tau),p(\tau),\lambda_0\big) d\tau \nonumber \\
               & & \label{SatzFreieZeitPMP1} - \sum_{j=1}^l \int_t^{t_{1*}} g_{jx}\big( \tau,x_*(\tau)\big) d\mu_j(\tau)
          \end{eqnarray}
          gen"ugt und die Transversalit"atsbedingung
          \begin{equation} \label{SatzFreieZeitPMP2} 
          p(t_{0*}) = {h_0}_x^T \big(t_{0*},x_*(t_{0*})\big) l_0.
          \end{equation}
          erf"ullt;
\item[b)] f"ur fast alle $t$ aus $[t_0,t_1]$ die Maximumbedingung
          \begin{equation} \label{SatzFreieZeitPMP3} 
          \chi_{\mathscr{A}_*}(t) \circ H\big(t,x_*(t),u_*(t),p(t),\lambda_0\big)
          = \max_{\substack{u^{i} \in U^{i}\\ 1\leq i\leq k}} H^{i}\big(t,x_*(t),u^{i},p(t),\lambda_0\big).
          \end{equation}
          gilt.
\item[c)] Die Hamilton-Funktion von beschr"ankter Variation und linksseitig stetig ist,
          und folgende Beziehungen gelten: 
          \begin{eqnarray}
          \lefteqn{\mathscr{H}\big(t,x_*(t),p(t),\lambda_0\big)
                      = \big\langle {h_1}_t\big(t_{1*},x_*(t_{1*})\big) , l_1 \big\rangle} \nonumber \\
          &-& \label{SatzFreieZeitPMP4} 
                \int_t^{t_{1*}} \chi_{\mathscr{A}_*}(\tau) \circ H_t\big(\tau,x_*(\tau),u_*(\tau),p(\tau),\lambda_0\big) d\tau 
              + \sum_{j=1}^l \int_t^{t_{1*}} g_{jt}\big( \tau,x_*(\tau)\big) d\mu_j(\tau), \\
          \lefteqn{\mathscr{H}\big(t_{0*},x_*(t_{0*}),p(t_{0*}),\lambda_0\big)
          = \label{SatzFreieZeitPMP5} - \big\langle {h_0}_t\big(t_{0*},x_*(t_{0*})\big) , l_0 \big\rangle.}
          \end{eqnarray}
\end{enumerate}
\end{theorem}

Ebenso wie in Aufgaben mit fester Zeit ist die adjungierte Funktion $p(\cdot)$ von beschr"ankter Variation und linksseitig stetig.

\newpage
\lhead[\thepage \, Gekoppelte Beh"alter]{Optimale Multiprozesse}
\rhead[Optimale Multiprozesse]{Gekoppelte Beh"alter \thepage}
\section{Zeitoptimale Steuerung gekoppelter Beh"alter}
\begin{figure}[h]
	\centering
	\fbox{\includegraphics[width=4cm]{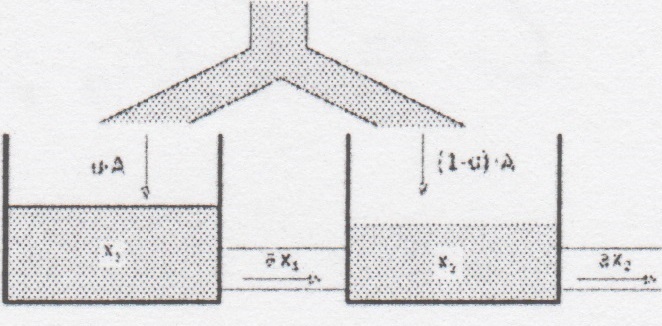}} \hspace*{2cm} \fbox{\includegraphics[width=4cm]{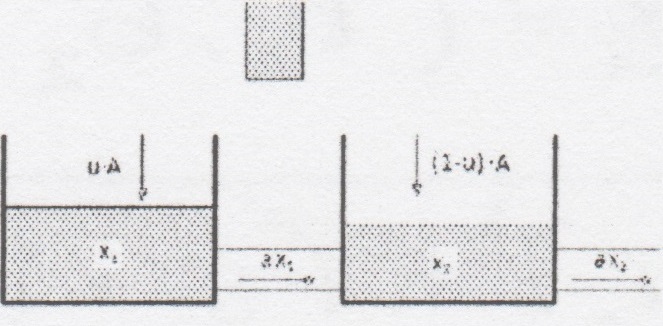}}
	\caption{Die Beh"alter in den verschiedenen Steuerungssystemen}
	\label{AbbildungContainer}
\end{figure}

Das System in Abbildung \ref{AbbildungContainer} zeigt zwei gekoppelte Beh"alter,
die mit einer Fl"ussigkeit gef"ullt sind.
In diesem Modell besteht die Kopplung der Beh"alter durch den zustandsabh"angigen Zuflu"s $ax_1$ vom ersten
in den zweiten Beh"alter.
Der zweite Beh"alter gibt seinen Inhalt mit der Rate $ax_2$ an die Umwelt ab.
Die beiden Steuerungssysteme bestehen einerseits aus dem Auff"ullen der Beh"alter mit konstanter Rate $A$,
die zwischen beiden Beh"altern durch die Steuerung $u^1$ aufgeteilt wird.
Andererseits kann die Auff"ullung komplett gestoppt werden.
Dies liefert folgende Dynamiken in den jeweligen Steuerungssystemen:
$$\begin{array}{l} \dot{x}_1 = -ax_1+u^1A, \\ \dot{x}_2=ax_1-ax_2+(1-u^1)A, \end{array}
  \qquad\qquad \begin{array}{l} \dot{x}_1=-ax_1,  \\ \dot{x}_2=ax_1-ax_2. \end{array}$$
Dabei beachte man, dass die ungesteuerte zweite Dynamik kein Spezialfall der ersten Dynamik ist.
In diesem Modell besteht das Ziel darin,
die F"ullst"ande $(x_1^0,x_2^0)$ in k"urzester Zeit in den Sollzustand $(x_1^f,x_2^f)$ zu "uberf"uhren.
Die Aufgabe des zeitoptimalen Multiprozesses lautet damit:
\begin{eqnarray}
&& \label{Container1} J\big(x(\cdot),u(\cdot),\mathscr{A} \big)
           = \int_0^T 1 \, dt \to \inf, \\
&& \label{Container2} \dot{x}(t) = \chi_{\mathscr{A}}(t) \circ \varphi \big(t,x(t),u(t) \big), \\
&& \label{Container3} \big(x_1(0),x_2(0)\big)=(x_1^0,x_2^0), \qquad \big(x_1(T),x_2(T)\big)=(x_1^f,x_2^f), \\
&& \label{Container4} u^1(t) \in [0,1], \quad \mathscr{A} \in \mathscr{Z}^2([0,T]), \quad a,A>0.
\end{eqnarray}
Wir wenden Theorem \ref{SatzFreieZeitPMP} an:
Die Pontrjagin-Funktion der Aufgabe (\ref{Container1})--(\ref{Container4}) lautet
$$\mathscr{A} \circ H(t,x,u,p,q,\lambda_0\big)=
  \left\{ \begin{array}{ll}
          p[-ax_1+u^1A]+q[ax_1-ax_2+(1-u^1)A]-\lambda_0, & t \in \mathscr{A}_1, \\[1mm]
          p[-ax_1]+q[ax_1-ax_2]-\lambda_0, & t \in \mathscr{A}_2. 
         \end{array}\right.$$
Unabh"angig von $\mathscr{A}$ erhalten wir daraus f"ur die Adjungierten $p(\cdot),q(\cdot)$:
\begin{equation} \label{Container5}
\bigg(\begin{array}{l} \dot{p}(t) \\ \dot{q}(t) \end{array}\bigg) =
  a \cdot \bigg(\begin{array}{r} p(t)-q(t) \\ q(t) \end{array}\bigg) \qquad\Rightarrow\qquad
  \bigg(\begin{array}{l} p(t) \\ q(t) \end{array}\bigg) =
  e^{at} \cdot \bigg(\begin{array}{l} p_0-atq_0 \\ q_0 \end{array}\bigg).
\end{equation}
Die Bedingung (\ref{SatzFreieZeitPMP3}) liefert Maximierungsaufgabe
\begin{equation} \label{Container6} \max_{i \in \{1,2\}} \Big\{ \max_{u^1\in [0,1]} [u^1p(t)+(1-u^1)q(t)],\; 0 \Big\}. \end{equation}
Au"serdem, da die Aufgabe autonom ist,
folgt $\mathscr{H}\big(t,x_*(t),p(t),\lambda_0\big) \equiv 0$ aus (\ref{SatzFreieZeitPMP4}) und (\ref{SatzFreieZeitPMP5}). \\[2mm]
Zur Auswertung der Maximierungsaufgabe (\ref{Container6}) sind die Relation zwischen $p(t),q(t)$ und $r(t)\equiv 0$ zu kl"aren.
Zun"achst gilt $(p_0,q_0)\not=(0,0)$ in (\ref{Container5}), denn sonst w"urde f"ur
\begin{enumerate}
\item[i)] $\lambda_0=0$ der Fall trivialer Multiplikatoren in Theorem \ref{SatzFreieZeitPMP} bzw.
\item[ii)] $\lambda_0\not=0$ ein Widerspruch zu $\mathscr{H}\big(t,x_*(t),p(t),\lambda_0\big) \equiv 0$
\end{enumerate}
auftreten.
Dann ergeben sich in Abh"angigkeit von $q_0 \in \R$ folgende Situation f"ur die Adjungierten,
wobei $\sigma$ den (m"oglichen) Schnittpunkt der Graphen von $p(\cdot)$ und $q(\cdot)$ bezeichnet:
\begin{enumerate}
\item[1)] $q_0=0$: $p(t)>q(t)\equiv 0$ f"ur $p_0>0$; $p(t)<q(t)\equiv 0$ f"ur $p_0<0$.
\item[2)] $q_0>0$: $p(t)<q(t)$ f"ur $p_0\leq q_0$ und $t>0$; es existiert $\sigma>0$ f"ur $p_0>q_0$.
\item[3)] $q_0<0$: $p(t)>q(t)$ f"ur $p_0\geq q_0$ und $t>0$; es existiert $\sigma>0$ f"ur $p_0<q_0$.
\end{enumerate}
Beachten wir nun weiterhin die optimale Stoppzeit $T_*$ und die (m"ogliche) Nullstelle $\sigma_0$ von $p(\cdot)$,
dann ergeben sich folgende Wechselstrategien $\mathscr{A}_*$ und Steuerungen $u^1_*(t)$:
\begin{enumerate}
\item[a)] Es ist nur ein Steuerungssystem aktiv auf $[0,T_*]$, wenn
          \begin{enumerate}
          \item[i)] $p(t)>\max\{q(t),0\}$ f"ur alle $t \in (0,T_*)$ und es ist
                    $\big(u^1_*(t),\chi_{\mathscr{A}^1_*}(t),\chi_{\mathscr{A}^2_*}(t)\big) \equiv (1,1,0)$;
          \item[ii)] $q(t)>\max\{p(t),0\}$ f"ur alle $t \in (0,T_*)$ und es ist
                     $\big(u^1_*(t),\chi_{\mathscr{A}^1_*}(t),\chi_{\mathscr{A}^2_*}(t)\big) \equiv (0,1,0)$;
          \item[iii)] $0>\max\{p(t),q(t)\}$ f"ur alle $t \in (0,T_*)$ und es ist
                      $\big(u^1_*(t),\chi_{\mathscr{A}^1_*}(t),\chi_{\mathscr{A}^2_*}(t)\big) \equiv (1,1,0)$.
          \end{enumerate}
\item[b)] Es gibt einen Strategiewechsel,
          \begin{enumerate}
          \item[i)] und zwar wenn $p_0>q_0>0$ und $\sigma \in (0,T_*)$:
                    $$\big(u^1_*(t),\chi_{\mathscr{A}^1_*}(t),\chi_{\mathscr{A}^2_*}(t)\big) \equiv (1,1,0)
                     \quad\longrightarrow \quad \big(u^1_*(t),\chi_{\mathscr{A}^1_*}(t),\chi_{\mathscr{A}^2_*}(t)\big) \equiv (0,1,0);$$
          \item[ii)] und zwar wenn $p_0<q_0<0$ und $\sigma_0 \in (0,T_*)$:
                     $$\big(\chi_{\mathscr{A}^1_*}(t),\chi_{\mathscr{A}^2_*}(t)\big) \equiv (0,1)
                      \quad\longrightarrow \quad\big(u^1_*(t),\chi_{\mathscr{A}^1_*}(t),\chi_{\mathscr{A}^2_*}(t)\big) \equiv (1,1,0).$$
          \end{enumerate}
\item[c)] Der Strategiewechsel ist nicht eindeutig bestimmbar, wenn $p_0<q_0=0$:
          $$\big(u^1_*(t),\chi_{\mathscr{A}^1_*}(t),\chi_{\mathscr{A}^2_*}(t)\big) \equiv (0,1,0)
            \quad\longleftrightarrow \quad  \big(\chi_{\mathscr{A}^1_*}(t),\chi_{\mathscr{A}^2_*}(t)\big) \equiv (0,1).$$
\end{enumerate}

\begin{figure}[h]
	\centering
	\fbox{\includegraphics[width=4cm]{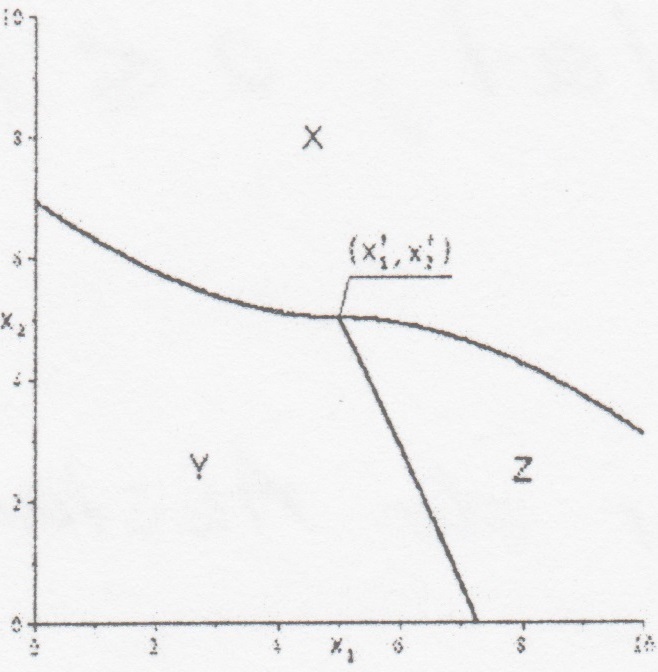}} \hspace*{2cm} \fbox{\includegraphics[width=4cm]{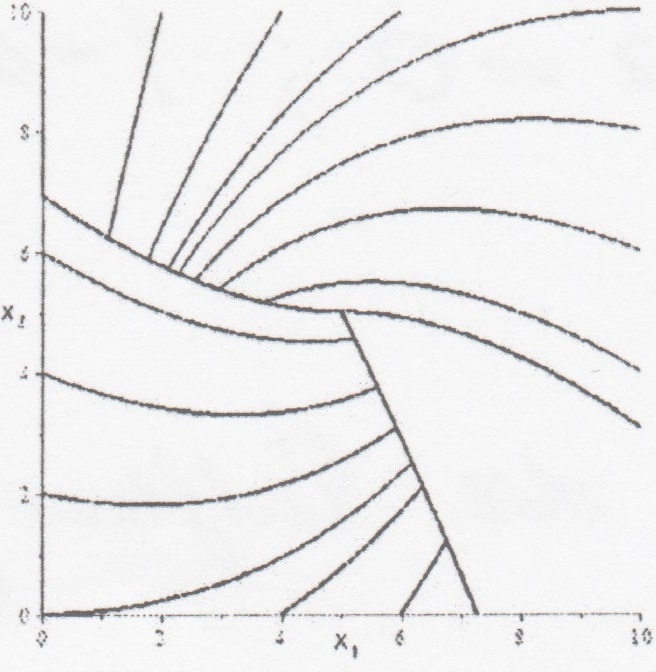}}
	\caption{Die Einteilung des Phasendiagramms und die optimalen Trajektorien}
	\label{AbbildungBereiche}
\end{figure}

Zur Illustration dieser "Ubersicht sei $a=1$, $A=10$, $x_1^0,x_2^0 \in [0,10]$ und $(x_1^f,x_2^f)=(5,5)$.
Die Zust"ande $\big(x_*^1(t),x_*^2(t)\big)$,
die zu den Steuerungen und aktivem System im Fall a) geh"oren,
teilen das Phasendiagramm $(x_1,x_2) \in [0,10] \times [0,10]$ in die Bereiche $X,Y,Z$ ein (Abbildung \ref{AbbildungBereiche}).

\begin{enumerate}
\item[$X$:] Die optimale Steuerung und Wechselstrategie wird durch ii) im Fall b) gegeben.
          Zum Zeitpunkt $t=\sigma_0$ erreicht die Trajektorie $\big(x_*^1(t),x_*^2(t)\big)$ den Rand des Bereiches $X$.
\item[$Y$:] Die optimale Steuerung und Wechselstrategie wird durch i) im Fall b) gegeben.
          Zum Zeitpunkt $t=\sigma$ erreicht die Trajektorie $\big(x_*^1(t),x_*^2(t)\big)$ den Rand des Bereiches $Y$.
\item[$Z$:] Die optimale Strategie ist nicht eindeutig, denn es tritt der Fall c) ein.
          Da in diesem Bereich die optimale Stoppzeit $T_*$ durch die schnellste Reduktion des Zustandes $x_1^0$ zu $x^f_1$
          festgelegt ist, muss im zweiten Beh"alter lediglich die Bedingung $x_2(T_*)=x_2^f$ sichergestellt werden.
\end{enumerate}

Wir f"ugen dem Beh"altermodell (\ref{Container1})--(\ref{Container4}) eine Zustandsbeschr"ankung hinzu:
$$x^2(t) \geq S \mbox{ f"ur alle } t \in [0,T], \qquad x_0^2,x_0^f,A>S.$$
Sei nun $\big(x_*^1(t),x_*^2(t)\big)$ eine optimale Trajektorie in der Aufgabe (\ref{Container1})--(\ref{Container4})
ohne Zustandsbeschr"ankung mit $x_*^2(t)<S$ f"ur ein $t \in (0,T_*)$.
Dann existieren Zeitpunkte $0<t_1<t_2<T_*$ mit $x^2(t_1)=x^2(t_2)=S$.
Nach Theorem \ref{SatzFreieZeitPMP} erhalten wir auf dem Teilst"uck $[t_1,t_2]$,
auf dem die Zustandsbeschr"ankung aktiv wird, f"ur die Adjungierten
$$p(t)=q(t)>0, \qquad \dot{p}(t)=\dot{q}(t)=0.$$
Dies liefert unter dem Ansatz $d\mu(t)=\lambda(t)\,dt$ f"ur das Borelsche Ma"s
$$\lambda(t)=a\big(2q(t)-p(t)\big)>0.$$
Eingeschr"ankt auf das Teilst"uck $[t_1,t_2]$ erhalten wir f"ur den Multiprozess
$$x^1_*(t)=x^1_*(t_1)+(t-t_1)(A-S), \quad x^2_*(t)=S, \quad u^1_*(t)= \frac{x^1_*(t)-S+A}{A}, \quad 
  \big(\chi_{\mathscr{A}^1_*}(t),\chi_{\mathscr{A}^2_*}(t)\big)= (1,0).$$
F"ur $a=1$, $A=10$, $x_1^0,x_2^0 \in [0,10]$, $(x_1^f,x_2^f)=(5,5)$ und $S=3,5$ ergibt sich das Phasenportrait
\begin{figure}[h]
	\centering
	\fbox{\includegraphics[width=4cm]{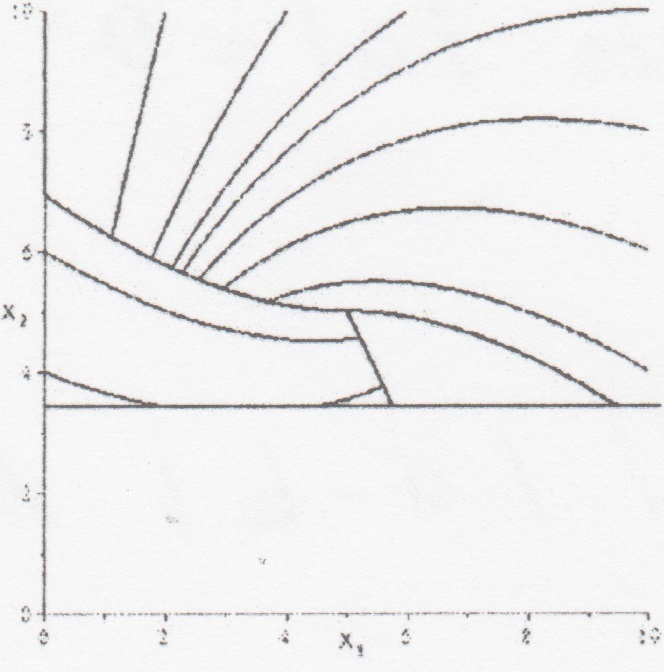}}
	\caption{Die optimalen Trajektorien im Phasendiagramm unter Zustandsbeschr"ankung}
	\label{AbbildungBereiche2}
\end{figure}
  
\newpage
\addcontentsline{toc}{section}{Literatur}
\lhead[\thepage \, Literatur]{Optimale Steuerung mit unendlichem Zeithorizont}
\rhead[Optimale Steuerung mit unendlichem Zeithorizont]{Literatur \thepage}

\end{document}